\documentclass[letterpaper,a4paper,oneside]{article}
\usepackage{graphicx} 
\usepackage{amsmath} 
\usepackage{amsxtra} 
\usepackage{amssymb} 
\usepackage{amsthm} 
\usepackage{latexsym} 
\usepackage{setspace} 
\usepackage[margin=1in]{geometry} 
\usepackage[titles]{tocloft} 
\usepackage{tikz}
\usepackage{tikz-cd}
\usepackage{verbatim}
\usepackage{tkz-graph}
\usepackage{adjustbox}
\usepackage{blindtext}
\usepackage{array}
\usepackage{fancyhdr,graphicx,amsmath,amssymb}
\usepackage{fancybox}
\usepackage{pseudocode}
\usetikzlibrary{positioning}
\usepackage{enumitem}
\usepackage{scrextend}
\usepackage{booktabs}
\usepackage{cite}
\usetikzlibrary{cd}
\usepackage{algorithmic}
\usepackage[T1]{fontenc}
\usepackage{hyperref}
\hypersetup{
    colorlinks=true,
    linkcolor=blue,
    filecolor=magenta,      
    urlcolor=blue,
    citecolor=black
}

\newcommand{\graph}{$G=(V,E)$ }

\usepackage[labelfont=bf,textfont=md]{caption}

\usepackage{titlesec}
\titleformat*{\section}{\large\bfseries}
\titleformat*{\subsection}{\bfseries}

\usepackage{authblk}
\title{\Large \textbf{\Large{Generalized Bijective Maps between $G$-Parking Functions, Spanning Trees, and the Tutte Polynomial}}}
\author{Carrie Frizzell}
\affil{Department of Mathematics, Kansas State University, 138 Cardwell Hall
    Manhattan, KS, 66506, USA; \href{mailto:cfrizz@ksu.edu}{cfrizz@ksu.edu}}
\date{}
\providecommand{\keywords}[1]
{
  \small	
  \textbf{\textit{Keywords---}} #1
}

\begin{document}
\maketitle
\noindent \textbf{Abstract.} We introduce an object called a tree growing sequence (TGS) in an effort to generalize bijective correspondences between $G$-parking functions, spanning trees, and the set of monomials in the Tutte polynomial of a graph $G$. A tree growing sequence determines an algorithm which can be applied to a single function, or to the set $\mathcal{P}_{G,q}$ of $G$-parking functions. When the latter is chosen, the algorithm uses splitting operations - inspired by the recursive defintion of the Tutte polynomial - to iteratively break $\mathcal{P}_{G,q}$ into disjoint subsets. This results in bijective maps $\tau$ and $\rho$ from $\mathcal{P}_{G,q}$ to the spanning trees of $G$ and Tutte monomials, respectively. We compare the TGS algorithm to Dhar's algorithm and the family described by Chebikin and Pylyavskyy in 2005. Finally, we compute a Tutte polynomial of a zonotopal tiling using analogous splitting operations.\medskip

\noindent\keywords{Tutte polynomial, G-parking functions, zonotopes}

\section{Introduction}
\label{intro}

To fix notation, given a multigraph $G=(V,E)$, label the vertices $V=\{q, v_1,\dots,v_n\}$, where $q$ 
is the root. The vertex and edge set will often be specified by $V(G),E(G)$ in context. If there are multiple edges between two vertices, order them. In each rooted subtree $T$ of $G$, we direct edges toward the root. When necessary, $e_h$ and $e_t$ are used for the head and tail of a directed edge $e=(e_h,e_t)$. Recall that a \textbf{spanning tree} of $G$ is a spanning, connected subgraph with $|V|-1$ edges.\medskip
\\
\textbf{Definition 1.1: }The \textbf{outdegree with respect to $A\subseteq V$}, denoted $outdeg_A(v)$, is the number of neighbors of $v$ not in $A\subseteq V$, with multiplicity.\smallskip\\
\textbf{Definition 1.2: }A \textbf{$G$-parking function} is a function $f:V(G)-\{q\}\rightarrow \mathbb{Z}_{\geq 0}$ such that any subset $A\subseteq V-\{q\}$ contains a vertex $v$ with $0\leq f(v)<outdeg_A(v)$.\medskip

\begin{figure}[htb]
    \centering
    \begin{tikzpicture}
    
    \draw (0,0)--(1,1);
    \draw (0,0)--(1,-1);
    \draw (1,1)--(1,-1);
    \draw (1,1)--(2,0);
    \draw (1,-1)--(2,0);
    \draw (2,0)--(3.5,0);
    
    \end{tikzpicture}
    
    \caption{The Tutte polynomial for the above graph is $T(G;x,y)=x^2+2x^3+x^4+xy+2x^2y+xy^2$. The multi-set $\mathcal{M}_G=\{x^2,x^3,x^3,x^4,xy,x^2y,x^2y,xy^2\}$.}
    \label{multi}
\end{figure}
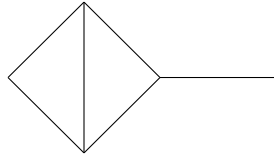

We write $f=(f(v_1),\dots,f(v_n))$. Let $\mathcal{P}_{G,q}$ denote the set of parking functions on $G$ with respect to $q$. 
Let $G-e$ mean deleting the edge $e$ from $G$. Contracting $G$ at $e$ means to delete $e$, then identify the endpoints of $e$. Denote contraction by $G/e$.\smallskip
\\
\textbf{Definition 1.3: }The \textbf{Tutte polynomial} $T(G;x,y)$ of G is the universal \textit{Tutte-Grothendieck graph isomorphism invariant} satisfying the following deletion/contraction principal, and defining $T(\bullet;x,y)=1$, $\bullet$ the graph with one vertex.
\begin{equation}T(G;x,y)=\begin{cases} 
      yT(G-e;x,y) & e \text{ a loop} \\
      xT(G/e;x,y) & e \text{ a bridge} \\
      T(G-e;x,y)+T(G/e;x,y) & \text{otherwise} 
   \end{cases}
\end{equation}
An equivalent definition is a closed formula over all spanning subgraphs of $G$. Let $c(A)$ be the number of connected components of a spanning subgraph $A$. Then
\begin{equation}T(G;x,y)=\sum_{A\subseteq G}(x-1)^{c(A)-c(G)} (y-1)^{|E(A)|+c(A)-|V|}\end{equation}
The symbol $\mathcal{M}_G$ will be used to indicate the multi-set with elements the terms of the Tutte polynomial of $G$. See Figure \ref{multi} for an example.

We first focus on the case of finite graphs due to the beautiful bijective correspondences between the terms of the Tutte polynomial, spanning trees, and $G$-parking functions. Among these is Dhar's burning algorithm \cite{DHAR90}; see also \cite{CHEBIKIN2005}, \cite{CHANG2010231}, \cite{bernardi:hal-00088479}, \cite{CORI200344}. The burning algorithm is applied to a graph labeled by a function $f: V(G)-\{q\}\rightarrow \mathbb{Z}_{\geq 0}$. Start a fire at $q$, and imagine it burns any edge it reaches. In order to burn through a vertex $v$, it must first burn through $z=f(v)$ edges which are incident to $v$. If the fire is able to burn through more than $z$ edges incident to $v$, it will burn through the vertex. All vertices burn if and only if $f$ is a $G$-parking function.

We will describe an algorithm which from a $G$-parking function simultaneously produces a spanning tree $T_f$ of $G$ and a monomial $x^{\alpha}y^{\beta}\in \mathcal{M}_T$, through the application of an object called a \textit{tree growing sequence} $\Sigma$. This results in two bijections $\tau:\mathcal{P}_{G,q}\rightarrow \mathcal{T}_G$, and $\rho :\mathcal{P}_{G,q} \rightarrow \mathcal{M}_T$. We prove the main theorem in section $\ref{split}$.\medskip
\\
\textbf{Theorem:} \textit{The maps $\rho$ and $\tau$ are bijective.\medskip}

The algorithm which achieves this is based on operations which simultaneously split each of the sets $\mathcal{P}_{G,q}$, $\mathcal{T}_G$, and $\mathcal{M}_G$ into two disjoint subsets. We show that these splittings are coherent in that they eventually force $1-1$ correspondences between the sets. As applications of the theorem, sections \ref{dhar}, \ref{cp}, and \ref{cmy} describe how Dhar's algorithm with a total edge order $O_E$ \cite{CORI200344}, \textit{proper sets of tree orders} $\{\Pi_G\}$\cite{CHEBIKIN2005}, and \textit{process orders} \cite{CHANG2010231}, respectively, can be fit into our definition. Let $\{O_E\}$ be the set of edge orders on $G$. We will define the maps in the diagram below and prove that it commutes. An auxilliary result is the association, via $\rho$, of a monomial to a $G$-parking function for the family of bijections in \cite{CHEBIKIN2005}. This is evidence that the $TGS$ algorithm can give a non-arbitrary bijection from $\mathcal{P}_{G,q}$ to $\mathcal{M}_G$ for algorithmic bijections between $\mathcal{P}_{G,q}$ and $\mathcal{T}_G$.

\begin{center}
\begin{tikzcd}
    & & \{\tau\} & \\
	\{O_E\} \arrow[urr, "D"] \arrow[r, "R" below]& \{\Sigma\} \arrow[ur, "F" left] \\
    & \{\Pi_G\} \arrow[uur, "\Phi" right] \arrow[u, "\Omega" left]
\end{tikzcd}
\end{center}

The Tutte polynomial is defined more generally for a matroid; see \cite{BRYOX92} for a thorough survey. In section \ref{zonotopes}, we compute a polynomial for a cubical zonotopal tiling using similar splitting operations to the TGS algorithm, and show that it is the Tutte polynomial of a specific matroid. In particular, if the vector configuration associated to the tiling is a cographical matroid, then the polynomial is the Tutte polynomial of the underlying finite graph, and we obtain bijections between tiles, $\mathcal{M}_G$, and $\mathcal{T}_G$. We conclude in section \ref{end} with a discussion that relates zonotopal tilings of cographical matroids to the bijective maps $\rho$ and $\tau$.
\vspace{3mm}

\noindent \textbf{Acknowledgements.} The author expresses gratitude to Ilia Zharkov for his role in the development and discussion of the ideas in this paper, and for patiently reading drafts. 

The author also thanks Kevin Long, Jeremy Martin, and McCabe Olsen for discussions related to this work during the 2019 Graduate Research Workshop in Combinatorics.\footnote{The 2019 Graduate Research Workshop in Combinatorics was supported in part by NSF grant \#1923238, NSA grant \#H98230-18-1-0017, a generous award from the Combinatorics Foundation, and Simons Foundation Collaboration Grants \#426971 (to M. Ferrara), \#316262 (to S. Hartke) and \#315347 (to J. Martin).}


\section{The Tree Growing Sequence}
\label{tgs}

\subsection{Definition and the Main Algorithm}
\label{deftgs}
We define the central object of this paper, the \textit{tree growing sequence}.\medskip
\\
\textbf{Definition:} Given a connected graph $G=(V,E)$ and the set $\mathcal{S}$ of all subsets of $E(G)$ containing $q$ as a vertex, and containing $\{q\}$ as a single element set, a \textbf{tree growing sequence (TGS)} is a collection of tuples $$\Sigma=\{(S,\sigma_S:~\mathcal{H}_S \rightarrow ~E(S))\}$$ where $S\in \mathcal{S}$, $\sigma_S$ is a function from the set $\mathcal{H}_S$ of subsets of $S$ containing $q$, $\sigma_S(T)\notin E(T)$, and the graph on $\sigma_S(T)\cup T$ is connected.\medskip

We will sometimes use a slight abuse of terminology, and instead call an element $S\in \mathcal{S}$ or $T\subset S$ a subgraph. Furthermore, the name ``tree growing sequence'' is used because the set $T$ will always be a tree in the application of our algorithm. 

Given a tree growing sequence $\Sigma$ and a function $f:V(G)-\{q\}\rightarrow \mathbb{Z}$, we apply the following algorithm to the tuple $(f, S, U, X, \alpha, \beta)$, where $U\subseteq V(G)$, $X\subseteq E(G)$, and $\alpha,\beta \in \mathbb{Z}_{\geq 0}$. The result will be a tree $T_f=(U,X)$ and a monomial $x^{\alpha}y^{\beta}$. Beginning with $S=E(G)$, $U=\{q\}$, $X=\emptyset$, $\alpha=0$, and $\beta=0$, the edge $\sigma_G(\{q\})=(e_h,q)$ is added to $X$ and $e_h$ added to $U$ if $f(e_h)=0$ and $e_h\neq q$ ($e$ not a loop). Furthermore, when $e$ is a bridge of $G$, then $\alpha$ is increased by one. If $e_h=q$, so that $e$ is a loop, delete it and increase $\beta$ by one. If $f(e_h)\geq 1$, the value of $f(e_h)$ is reduced by one. The edge $e$ is not added to $X$, and we equate this with edge deletion by replacing $E(G)$ with $E(G)-e$. If it is the case that $f(e_h)<0$, we terminate the algorithm. 

In subsequent steps, we consider the tuple $(f, S, U, X, \alpha, \beta)$, where the value of $f$ at some vertices may have been reduced in previous steps. For each image $\sigma_S(T)=e$, we assume that $e_t$ is a vertex of $T$. If $e_h$ is also a vertex of $T$, then $e$ will be called a loop. The set $S=E(G)-\{e\}_D$, where the edges $\{e\}_D$ have been deleted. 
The algorithm is shown in below.

\vspace{1.5mm}


\begin{pseudocode}[doublebox]{Tree Growing Sequence Algorithm }{f,S,U,X,\alpha,\beta}
\textbf{Input:} \text{ A graph } G=(V,E) \text{ with root vertex } q, \text{ tree growing sequence } \Sigma,\\\text{ and an integer valued function } f \text{ on the vertices.}\\
\textbf{Output: } \text{A tree }T_f \text{ and monomial } x^{\alpha}y^{\beta}.\\
\textbf{Initialization:}\\
 S=G\\
 U=\{q\}\\
 X=\emptyset\\
 \alpha=0\\
 \beta=0\\
 T=(\{q\},\emptyset)\\

\WHILE \sigma_S(T) \text{ is defined }\DO\\
\BEGIN
    \IF f(\sigma_S(T))<0
        \THEN
        terminate
    \ELSEIF f(\sigma_S(T))=0 
    \THEN
    \BEGIN
        \IF e=\sigma_S(T) \text{ a bridge of }S
        \THEN
        \BEGIN
            \alpha \GETS \alpha+1\\
            X\GETS X\cup e\\
            U\GETS e_h\\
            S\GETS S\cup e\\
            T=(U,X)
        \END
        \ELSE
        \BEGIN
            X\GETS X\cup e\\
            U\GETS e_h\\
            S\GETS S\cup e\\
            T=(U,X)
        \END
    \END
\ELSEIF T\cup e \text{ not a tree}
\THEN
\BEGIN
\beta \GETS \beta+1\\
S \GETS S-e\\
\END
\ELSE
\BEGIN 
f(e) \GETS f(e)-1\\
S\GETS S-e\\
\END
\END\\
\OUTPUT{T_f=(U,X),\,x^{\alpha}y^{\beta}}
\end{pseudocode}


We illustrate in Figure \ref{dc_action} the possibilities for updating the tuple when applying the algorithm.
\vspace{2mm}

\begin{figure}[htb]
    \centering
    \begin{tikzcd}
        (f,S,U,X,\alpha,\beta) \arrow[rr,"a"]& &(f,S,U\cup e_h,X\cup e,\alpha,\beta)\\
        (f,S,U,X,\alpha,\beta) \arrow[rr,"b"] & & (f, S, U\cup e_h, X\cup e, \alpha+1,\beta)\\
        (f,S,U,X,\alpha,\beta) \arrow[rr,"c"]&& (f, S-e, U,X,\alpha,\beta+1)\\
        (f,S,U,X,\alpha,\beta) \arrow[rr, "d"]&& (f,S-e,U,X,\alpha,\beta)
    \end{tikzcd}
    \caption{Possibilities for updating the tuple.}
    \label{dc_action}
\end{figure}
\begin{figure}[htb]
    \centering
    \includegraphics[scale=0.53]{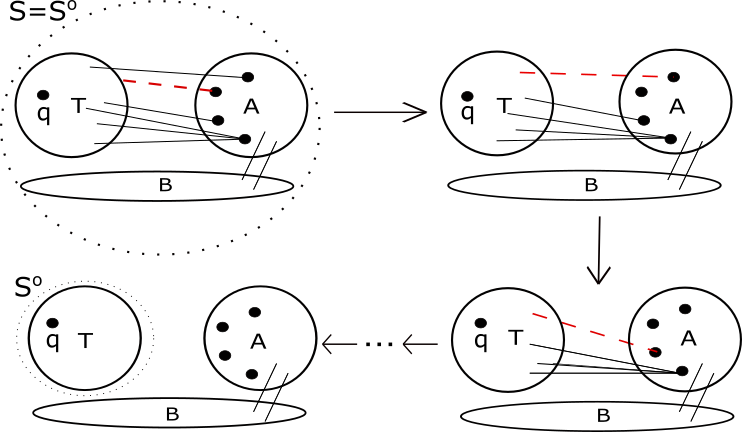}
    \vspace{1mm}
    
    \caption{The vertex set $B=V(G)-V(T)-A$. The picture shows what happens when $\sigma_S(T)=e$ for an edge joining a vertex in $A$ to one in $V(T)$.}
    \label{nonpark}
\end{figure}

\newpage
\noindent \textbf{Proposition 2.1.1:} \textit{For any tree growing sequence $\Sigma$, applying the $TGS$ algorithm to a function $f:V(G)-\{q\}\rightarrow \mathbb{Z}$ will terminate on a spanning tree of $G$ if and only if $f\in \mathcal{P}_{G,q}$.}
\begin{proof}
Fix a root vertex $q$ and $f\in \mathcal{P}_{G,q}$. If the algorithm terminates at non-spanning $T_f$, then $T_f$ spans $S$ but not $G$. This implies that $V(S)\neq V(G)$ and we have deleted all edges between $V(S)$ and $U=V(G)-V(S)$. Then we can find some $A\subseteq U$ such that $outdeg_A(v)\leq f(v)$ for all $v\in A$. However, this is impossible since $f\in \mathcal{P}_{G,q}$. Hence, $V(S)=V(G)$, and $T$ spans ~$G$.
 
Conversely, if $h\notin \mathcal{P}_{G,q}$, then a tree growing sequence will not terminate on a spanning tree of $G$. Let $A\subseteq V-\{q\}$ be a subset such that all vertices $v\in A$ satisfy $outdeg_A(v)\leq h(v)$. It will suffice to let $A$ consist of a single vertex $v$, because any such subset $A$ can be thought of as a single vertex with $deg(A)=\sum_{v\in A}outdeg_A(v)$. This translates to $0<deg(v)\leq h(v)$ (excluding loops). Consider the first time that $\sigma_S(T)=(v,u)$, $u\in V(T)$. This will eventually occur, because $\sigma_S(T)$ is defined as long as $T\neq S$. The edge $(v,u)$ will be deleted because it was assumed that $deg(v)>0$. Moreover, we reduce $h(v)$ by one. Every time $\sigma_{S'}(T')=(v,u')$, the edge will be deleted, and $h(v)$ reduced by one. Since $deg(v)\leq h(v)$, we will eventually exhaust all edges from $A$ to $T$. Hence, we will not get a spanning tree by applying $\Sigma$ to $h$ (Figure \ref{nonpark}). 
\end{proof}
\noindent We define the map $\tau:\mathcal{P}_{G,q}\rightarrow \mathcal{T}_G$ to be the assignments $f\mapsto T_f$ according to $\Sigma$.\medskip
\\
\textbf{Proposition 2.1.2: }\textit{If $f\in \mathcal{P}_{G,q}$, then the above algorithm always produces a monomial $x^{\alpha}y^{\beta}$ in the multiset $\mathcal{M}_T$ when applied to $f$.}
\begin{proof}
Start with $\alpha=\beta=0$. If $e=\sigma_S(T)$ is a bridge of $S$, then increase $\alpha$ by one. If $e$ is a loop, meaning $f(e_h)=0$ and $T\cup e$ has a cycle, then delete it, and increase $\beta$ by one. In light of equation $(1)$, we are simply isolating a monomial of $T(G;x,y)$ when computing it via recursion, and this is the monomial which we associate to $f$.
\end{proof}

\noindent The above proposition is nothing new. However, it is the starting point for a closed formula for the Tutte polynomial as a sum over $\mathcal{P}_{G,q}$ - done by Chang, Ma, and Yeh in \cite{CHANG2010231} - and serves as inspiration to generalize known algorithmic bijections. Also note that the set $X\subset E(S)$ can be viewed as contracted edges, though technically we do not alter the structure of the subgraph when adding an edge to $X$. 


\subsection{The Splitting of $\mathcal{P}_{G,q}$}
\label{split}

\begin{figure}[htb]
\centering
\adjustbox{scale=0.85,center}{
\begin{tikzcd}
& & &\arrow[dll] \sigma_G \arrow[drr] & & &\\
& \sigma_{G} \arrow[dl] \arrow[dr] & & & & \sigma_{G-e} \arrow[dl] \arrow[dr] & \\
\sigma_{G}& & \sigma_{G-e'} & & \sigma_{G-e}& &\sigma_{(G-e)-e'} \\
\updownarrow & \vdots & \vdots & \vdots & \vdots & \vdots & \updownarrow\\
\sigma_G & & & & & & \sigma_{((G-e)-e')-l}
\end{tikzcd}
}
\adjustbox{scale=0.85,center}{
\begin{tikzcd}
& & &\arrow[dll] \mathcal{P}_G \arrow[drr] & & &\\
& \mathcal{P}_{G/e} \arrow[dl] \arrow[dr] & & & & \mathcal{P}_{G-e} \arrow[dl] \arrow[dr] & \\
\mathcal{P}_{(G/e)/e'}& & \mathcal{P}_{(G/e)-e'} & & \mathcal{P}_{(G-e)/e'}& &\mathcal{P}_{(G-e)-e'} \\
\updownarrow & \vdots & \vdots & \vdots & \vdots & \vdots & \updownarrow\\
\mathcal{P}_{((G/e)/e')/b} & & & & & & \mathcal{P}_{((G-e)-e')-l}
\end{tikzcd}}
\caption{Binary trees illustrating how splitting the parking functions corresponds to the application of $\Sigma$. Here, $b$ means bridge and $l$ means loop. Note that the edges denoted $e'$ are not necessarily the same on each side of the tree. We can replace $\mathcal{P}$ with $\mathcal{M}$ and the splitting looks the same.}
    \label{bin}
\end{figure}
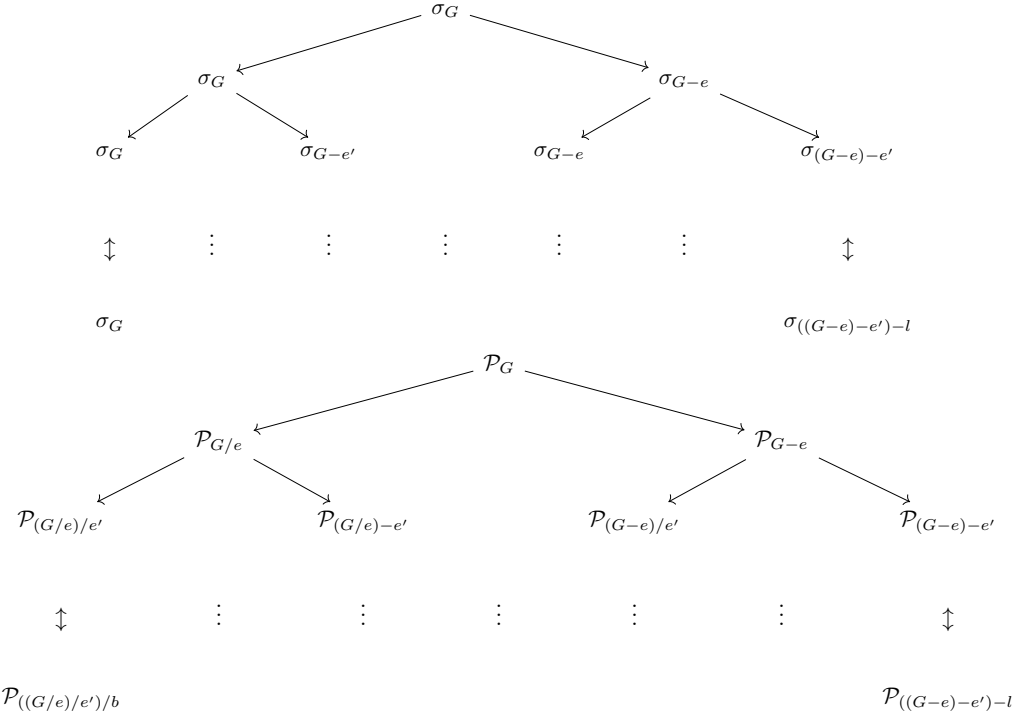

We change our philosophy from the previous section: instead of taking a single $G$-parking function $f$ and producing a spanning tree and monomial, we begin with the 
set of parking functions $P_{G,q}$ and perform splitting operations with respect to the deletion/contraction principle. That is, split the parking 
functions according to whether the edge $e=(e_h,e_t)$ is added to contracted (added to $X$) or deleted; see Figure \ref{bin} for a visual. This splitting will also result in the bijections $\tau:\mathcal{P}_{G,q}\rightarrow \mathcal{T}_G$ and $\rho :\mathcal{P}_{G,q} \rightarrow \mathcal{M}_G$. To this end, we include the proofs of two lemmas. We use the convention that when an edge $e=(e_h,e_t)$ is contracted and $e_h,e_t$ are identified, the ``thickened'' vertex is called $e_t$. We begin with letting $e=(e_h,q)=\sigma_G(\{q\})$ for an arbitrary TGS $\Sigma=\{(S,\sigma_S)\}$. 
\smallskip
\\
\textbf{Lemma 2.2.1:} \textit{If $e$ is a loop, then $\mathcal{P}_{G,q}=\mathcal{P}_{G-e,q}$.\medskip \\}
\textbf{Lemma 2.2.2:} \textit{If $e$ is a bridge, then $P_{G,q}$ is in one-to-one correspondence with $P_{G/e,q}$.\medskip \\}
\textbf{Lemma 2.2.3} \cite{CHANG2010231}: \textit{There is a bijection $\phi$ between the set of $G$-parking functions $f$ with $f(e_h)=0$ and the set of $(G/e)$-parking functions.\medskip \\}



\noindent \textbf{Lemma 2.2.4} \cite{CHANG2010231}: \textit{There is a bijection $\psi$ between the set of $G$-parking functions $f$ with $f(e_h)\geq1$ and the set of $(G-e)$-parking functions.\medskip \\}
\noindent \textbf{Corollary 2.2.3:} \textit{For any graph $S$ with fixed root $q$, there is a bijection between $\mathcal{P}_{S,q}$ and $\mathcal{P}_{S/e,q}\sqcup \mathcal{P}_{S-e,q}$.
\medskip\\}
Recall that the map $\tau:\mathcal{P}_{G,q}\rightarrow\mathcal{T}_G$ is the assignment of each $G$-parking function $f$ to the spanning tree $T_f$ on which a tree growing sequence $\Sigma$ terminates, and let $\rho :\mathcal{P}_{G,q} \rightarrow \mathcal{M}_T$ be the assignment of a monomial to each $f$. We will, in general, get different $\rho, \tau$ for different $\Sigma$.\medskip
\\
\textbf{Theorem 2.2.4:} \textit{The maps $\rho$ and $\tau$ are bijective.}
\begin{proof}
It is a well-known fact that the sizes of the three sets $\mathcal{P}_{G,q}, \mathcal{T}_G, \mathcal{M}_T$ are equal. Hence, it is enough to show that if $f,g\in \mathcal{P}_{G,q}$ are not equal, then $\tau(f)\neq \tau(h)$, and for each $x^{\alpha}y^{\beta}\in \mathcal{M}_T$, there is a unique (up to permuting identical elements) $f$ with $\rho(f)=x^{\alpha}y^{\beta}$.




Fix $f\neq h$. By Corollary 2.2.3, each splitting produces a bijection between $\mathcal{P}_{S,q}$ and $\mathcal{P}_{S/e,q}\sqcup \mathcal{P}_{S-e,q}$, where $S=G-\{e\}_D$ according to the edges previously deleted. As we never contract edges, we view $\mathcal{P}_{S/e,q}$ as the set of parking functions such that $e=\sigma_S(T)$ is added to $X$. If $\tau(f)=\tau(h)$, then the same set of edges $\{e^1,\dots,e^m\}$ is contracted in the paths for both. However, this implies that either $f$, $h$ have the same path, which implies $f=h$; o,r $f$ and $h$ split and have the same edges contracted. This is impossible, as there is some $e$ for which $f$ is in the contraction set, and $h$ is in the deletion set. Therefore, $e\in T_f$, but $e\notin T_h$, and $\tau(f)\neq \tau(h)$.

The statement that each monomial $x^{\alpha}y^{\beta}$ in the multiset $\mathcal{M}_T$ has a unique preimage $\rho^{-1}(x^{\alpha}y^{\beta})\in \mathcal{P}_{G,q}$ up to permuting repeated elements can be proven by splitting $\mathcal{M}_G$ using formula $(1)$ in section \ref{deftgs}. If $e$ is neither a bridge nor loop, then $\mathcal{M}_G=\mathcal{M}_{G/e}\sqcup \mathcal{M}_{G-e}$; if $e$ is a loop, then $\mathcal{M}_G=y\cdot \mathcal{M}_{G-e}$; and if $e$ is a bridge, $\mathcal{M}_G=x\cdot \mathcal{M}_{G/e}$. Hence, if $e=\sigma_S(T)$ is a loop or bridge, no splitting occurs. If $e$ is neither, then $\mathcal{M}_{G/e}$ corresponds to $\mathcal{P}_{G/e}$, and $\mathcal{M}_{G-e}=\mathcal{P}_{G-e}$. The result of iterating the process until it terminates is that each $f\in \mathcal{P}_{G,q}$ is in corresondence with a unique element of $\mathcal{M}_G$ (again, up to permuting identical monomials).
\end{proof}
\section{Comparison of TGS to Known Algorithms}
\label{frame}
This section is dedicated to relating tree growing sequences to formerly established bijective algorithms between the three objects of interest. We focus on Dhar's algorithm in Section \ref{dhar} and the family of bijections described by Chebikin and Pylyavskyy \cite{CHEBIKIN2005} in Section \ref{cp}. The bijection between $G$-parking functions and $\mathcal{M}_T$ given by Chang, Ma, and Yeh \cite{CHANG2010231} is discussed in Section \ref{cmy}.
\subsection{Global Edge Order and Dhar's Algorithm}
\label{dhar}
Given a global edge order $O_E:E(G)\rightarrow \{1,\dots, |E(G)|\}$, we construct a tree-growing sequence $\Sigma_{O_E}$ by defining for all $S\subseteq E(G)$ and subtrees $T\subseteq S$ the image $\sigma_S(T)=e$ to be the largest available 
edge which maintains a connected graph at each step. Call this construction the map $R:\{O_E\}\rightarrow \{\Sigma\}$ from the set of edge orders to the collection of tree growing sequences. This definition of $\Sigma_{O_E}$ mimics Dhar's burning algorithm ``with memory''; see, i.e. \cite{BAKER2013164}, \cite{CORI200344} for explicit algorithms and proofs of the Dhar bijection between $G$-parking functions (also referred to as $q$-reduced divisors) and spanning trees using a total edge order. In the notation for the TGS algorithm, Dhar's algorithm chooses the edge $e=max_{O_E}\{(v,u)\,|\,u\in U, v\notin U\}$. The edge $e$ is added to $X$ if $f(v)=|\{(v,u)\in E(G)-E(S)\}|$. Thus, the definition of $\Sigma_{O_E}$ is almost the same, except it may attempt to grow an edge which creates a cycle. 
Denote $D_{O_E}(f)$ the image of $f$ under Dhar's algorithm with edge order $O_E$. For a chosen root $q$, let $\Sigma_q$ denote the above TGS where we start at the root.\medskip
\\
\textbf{Proposition 3.1.1:} \textit{The map $R:\mathcal{O}_E\rightarrow \Sigma_q$ commutes with Dhar's algorithm, for any root ~$q$.} 
\begin{proof}

The map $R$ is defined as above. Fix a root $q$. If $q\in T\subseteq S$, then $\sigma_S(T)=e$, where $e=~max_{O_E}(E(S)-E(T))$ and $T\cup e$ is connected. If $e_h\notin U, e_t\in U$, then the edge is the same one chosen in Dhar's algorithm. Furthermore, the edge $e$ is deleted if $f(v)\geq 1$ (including after being reduced) which is equivalent to $f(v)>|\{e\in E(S)-X\,|\, v\in e\}|$. The edge is added to $X$ if precisely $f(v)$ edges incident to $v$ have been deleted. On the other hand, if $e_h, e_t\in U$, then $e$ will be deleted. Thus, we do not add this edge to $X$, and since such an edge is never considered in Dhar's algorithm - it is ignored - the diagram commutes.
\end{proof}

\begin{center}
\begin{tikzcd}
    \{O_E\} \arrow[r, "R"] \arrow[dr, "D" below left]& \{\Sigma\} \arrow[d, "F"]\\
    & \{\tau\},
\end{tikzcd}
\end{center}

\noindent Applying Dhar's algorithm to a $G$-parking function will also give a bijection with $\mathcal{M}_T$ via the notions of \textit{internal} and \textit{external activity} of the edges of $D_{\mathcal{O}_E}(f)$; this is how how Tutte originally defined the polynomial in \cite{tutte_1954}. An edge $e$ is internally active if it is smallest, according to $O_E$, in the unique cocircuit (cut-set) of $(G-T)\cup e$. Dually, an edge is externally active if it is the smallest in the unique cycle of $T\cup e'$. The Tutte polynomial can be written as a sum over $\mathcal{T}_G$: $$T(G;x,y)=\sum_{\mathcal{T}_G}x^{ia}y^{ea}$$ with \textit{ia} and \textit{ea} denoting the number of internal and external edges, respectively, of the tree $T$ according to $\mathcal{O}_E$. Commutativity of the diagram implies that $\tau(f)$ is the same monomial corresponding to $T_f$ in the above sum. On the other hand, one can ask if the internally and externally active edges match with the bridges and loops in the tree growing sequence algorithm.\medskip
\\
\textbf{Proposition 3.1.2:} \textit{If an edge $e\in E(G)$ is internally active for the tree $T_f$, then it is a bridge when added to $X$ during application of $\Sigma_{O_E}$.}
\begin{proof}
Say $e$ contributes to  the exponent $\alpha$, where $\tau(f)=x^{\alpha}y^{\beta}$. Then at some step of applying $\Sigma_{O_E}$ to $f$, $\sigma_S(T)=e$ is a bridge of $S$. Hence, either $e$ is a bridge of $G$, or there is a circuit $C$ in $S$ of which $e$ is the smallest among any adjacent edge $e'$ in $G-S$ - i.e. edges which have already been deleted in the tree growing process. Then it is the smallest edge in the unique cocircuit $B$ of $(G-T_f)\cup e$ containing $e$ and any $e'$ as described above.
\end{proof}

\noindent \textbf{Corollary 3.1.3: }\textit{The following diagram commutes.}
\begin{center}
\begin{tikzcd}
    \{O_E\} \arrow[r, "R"] \arrow[dr, "K" below left]& \{\Sigma\} \arrow[d, "F"]\\
    & \{\rho\}
\end{tikzcd}
\end{center}

\subsection{Proper Sets of Tree Orders}
\label{cp}

In \cite{CHEBIKIN2005}, a family of bijections between $G$-parking functions and spanning trees is produced using an object called a \textit{proper set of tree orders}, $\Pi_G$. Let \graph be a graph and choose a labeling of the vertices $\{v_1,\dots,v_n\}$. Given an ordering $\pi(T)$ on the vertices of every sutree $T$ rooted at $q$, the collection $\Pi_G=\{\pi(T)\:|\: T\subset G \text{ a rooted tree}\}$ is a proper set of tree orders if and only if 
the orders are compatible in the obvious way on overlaps (rooted at $q$) and a directed edge $(u,v)\in T$ means $v<u$ in $\pi(T)$. Specifically, the former translates to if the overlap of $T$ and $T'$ contains a rooted tree, and $i,j$ are vertices in this overlap, then $i<_{\pi(T)}j\iff i<_{\pi(T')}j$. Let $\pi(T)(q)=0$ for any $T$. Note that if the trees $T$ and $T'$ differ only by a choice of a set of multi-edges, the orders $\pi(T)$ and $\pi(T')$ must be the same. Examples include 
tree orders induced by vertex orderings constructed by \textit{breadth-first}, \textit{depth-first}, and \textit{vertex adding} algorithms. These three orders can all be constructed from the example below.
\begin{description}
\item \textbf{Example 3.2.1: \cite{CHEBIKIN2005}} One way to construct $\Pi_G$ is from a partial order on the set of (open) paths ending at $q$. The partial order must satisfy the conditions
that paths which intersect along another path at $q$ are comparable, and $A\preceq A\cup <v_k,\dots,v_k'>$.
\item The partial order $\preceq$ descends to a proper set of tree orders $\Pi_{\preceq}$. Given any rooted subtree $t\subset G$, and
distinct vertices $v, w\in V(T)$, the order $\pi_{\preceq}(T)$ is determined by the ordering of the paths from $q$ to $v$ and from $q$ to $w$.
Since these paths intersect along a path starting at $q$, they are comparable. However, not every $\Pi_G$ arises in this manner, 
see \cite{CHEBIKIN2005}.
\end{description}

\noindent We define a map $\Omega: \{\Pi_G\}\rightarrow \{\Sigma\}$. Fix $\Pi_G$. Consider any subgraph $S\in \mathcal{S}$, and any rooted subtree $T\subseteq S$. Then define $\sigma_S(T)=e$ according to the following:
\begin{enumerate}[label=(\roman*)]
    \item \begin{enumerate}
        \item Take the smallest edge according to $\pi(T)$ from every vertex a neighbor of $T$. Call this tree $T'$.
        \item Let $\sigma_S(T)=e$ be the edge in $T'$ such that $e_h$ is the smallest vertex in $V(T')-V(T)$ according to $\pi(T')$.
    \end{enumerate}
    \item If there is no edge in $S$ which satisfies $(i)$, let $\sigma_S(T)=e'$ for the smallest possible edge $e'$ induced by $\pi(T)$ such that $T\cup e'$ is connected.
\end{enumerate}

If no edge satisfies $(i)$ or $(ii)$, $\sigma_S(T)$ is undefined. Again, this happens when $T$ is equal to the connected component of $S$ containing $q$.

For example, given $\pi(T)=\{q,\dots,u_r\}$ where $u_j$ is the $j-th$ vertex in the order, defining $\sigma_T$ as above ensures that we grow $T$ according to the order $\pi(T)$. That is, $$\sigma_T(T_k)=e, T_k=(V_k=\{q,u_1,\dots,u_k\}, E_k), e=(u_{k+1},u)\in T, u\in V_k\}.
$$

\noindent \textbf{Proposition 3.2.1:} \textit{The map $\Omega$ is an injection from the collection of proper sets of tree orders to the collection of tree 
growing sequences.}

\begin{proof}

For each $T\subseteq S\subseteq G$, there is a unique image $\sigma_S(T)=e$, when defined. If not, there are two edges $e,e'$ satisfying the conditions. This means $e, e'$ are both minimal according to either $(i)$ or $(ii)$, which is impossible. Assembling this data into maps $\sigma_S$ and letting $\sigma_H$ be undefined for $q\notin H\subset G$ is precisely the data of a tree growing sequence $\Sigma=\{(S,\sigma_S)\}$.\smallskip

To show injectivity, we must show that if $\Omega(\Pi^a_G)=\Omega(\Pi^b_G)$, then $\Pi_G^a = \Pi_G^b$. Suppose otherwise. Then there is a rooted 
subtree $T'\in G$ such that $\pi_G^a(T')\neq \pi_G^b(T')$. Assume that $\Pi_G^a$ and $\Pi_G^b$ differ at the $k$-th vertex, i.e. $u_k^a\neq u_k^b$. 
Then $\sigma^a_{T'}(T'_{k-1})\neq \sigma^b_{T'}(T'_{k-1})$, which implies that $\Omega(\Pi_G^a)\neq \Omega(\Pi_G^b)$. 
Therefore, $\Omega$ is injective.  
\end{proof}
\noindent \textbf{Example 3.2.2:} We will borrow an example from \cite{CHEBIKIN2005}, pp 33-34, where $\Pi_G$ is the proper set of tree orders such that $i<_{\pi(t)} j$ if either $d_t(q,i)<d_t(q,j)$, or the distances are equal and $i<j$ in $G$.  Several cases are presented.
\begin{description} 
\item In case 1, we have the subtree $t$ of $T_1$ (left) and $T_2$ (right) shown with dotted edges. If $S=G$ with vertex order given, then we must have $\sigma_G(t)=(2,1)$. If we delete $(2,1)$, we have the subgraph $S$ (below $G$), and $\sigma_S(t)=(2,3)$.
\item In case 2, consider the subtree $t'$. We need to know how to grow $t'$ - if at all - in a given subgraph. First, let $S=G$. All spanning trees containing $t'$ are shown. 
We can check that we must define $\sigma_G(t')=(2,1)$. If we remove $(2,1)$, the map $\sigma_S, S=G-(2,1)$ will have image $(3,1)$ when applied to $t'$. The graphs to the right of $S$ are maximal subtrees of $S$.
\item One may observe that in light of the definition of $\sigma_G(T_2)$, the definition of $\sigma_G(T_1)$ is excess data, because we would be deleting the edge $(2,1)$ before growing the edge $(2,3)$. However, we want to define $\sigma_S$ on all edge subsets $S$ which form a connected subgraph, whether or not the data will be needed when applying $\Omega(\Pi)$.
\end{description}
\begin{figure}[htb]
    \centering
   \includegraphics[scale=0.45]{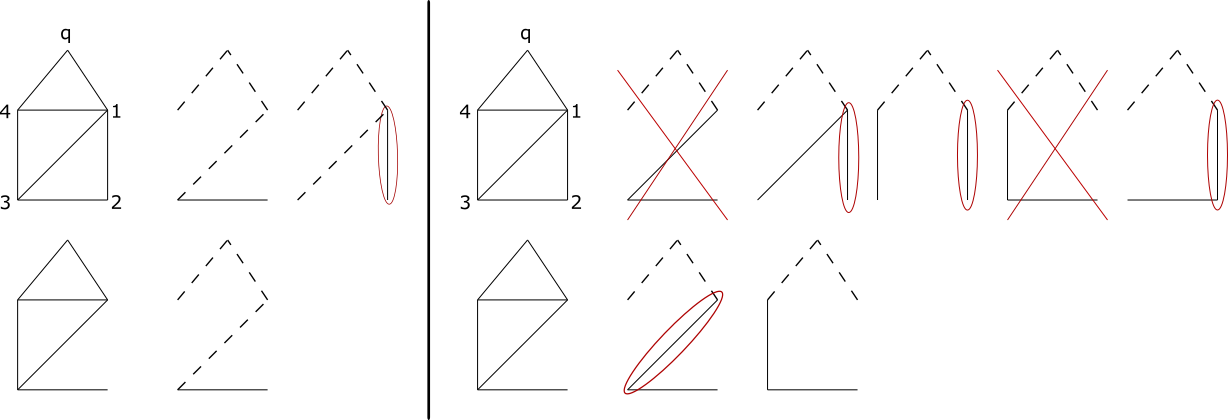}
    \caption{Some elements of the tree growing sequence $\Omega(\Pi_G)$.}
    \label{example1}
\end{figure}
\noindent The above proposition establishes that any proper set of tree orders can tell us how to proceed with the TGS algorithm.
However, it is desirable to have commutativity of the diagram in the theorem below. Before the theorem, we describe the bijective map $\Phi_{\Pi}: \mathcal{P}_{G,q}\rightarrow \mathcal{T}_G$, first given in \cite{CHEBIKIN2005}. Fix $f\in \mathcal{P}_{G,q}$. Declare $p_0=q$ and $T_0=\{q\}$. At each step $k$, let $T_{k-1}$ be the current subtree grown. The next edge to be grown, denoted $e_k=(p_k,v),v\in V(T_k)$, is the one that satisfies these conditions:
\begin{enumerate}
    \item There are at least $f(p_k)+1$ edges from $p_k$ to $T_{k-1}$,
    \item The edge $e_k$ is larger than precisely $f(p_k)$ of these edges, and
    \item The vertex $p_k$ is minimal among all vertices with edges satisfying $(i),(ii)$, according to the order of the tree obtained from $T_{k-1}$ by adjoining these edges.
\end{enumerate}

\noindent The labels $p_0,\dots,p_n$ comprise the order $\pi(T_f)$, in that $p_0<_{\pi(T)}\dots<_{\pi(T)}p_n$ (\cite{CHEBIKIN2005}, \textit{Lemma 2.3}).\medskip
\\
\textbf{Theorem 3.2.2: }\textit{The following diagram commutes.}
\begin{center}

\begin{tikzcd}
    \{\Pi_G\} \arrow[r, "\Omega" ] \arrow[dr, "\Phi" below left]& \{\Sigma\} \arrow[d, "F"]\\
    & \{\tau\}
\end{tikzcd}
\end{center}

\begin{proof}

Fix $f\in \mathcal{P}_{G.q}$ and $\Pi\in \{\Pi_G\}$. It will be shown that $\Omega(\Pi)(f)=\Phi_{\Pi}(f)$. We will argue that if an edge is added to $X$ when applying $\Omega(\Pi)$ to $f$, then it is in $\Phi_{\Pi}(f)$. Since we know $T_f=(V,X)$ is spanning by Proposition 2.1.1, this will prove the claim.

Consider the algorithm for constructing $\Phi_{\Pi}(f)=T_n$. At step $k$, let $V_k$ be the vertices not in $T_{k-1}$, $U_k\subseteq V_k$ the vertices adjacent to some vertex in $T_{k-1}$, and $W_k$ the set of vertices satisfying $(1)$. For $k=1$, we consider vertices with at least $f(v)+1$ edges to $q$. The edge $(v,q)$ satisfying condition $(2)$ will be in $E(T_n)$, for all $v\in W_1$. This is because $\pi(T)(q)=0$ for all $T$, so any edge from $v$ to future vertices in $T$ is larger than $(v,q)$. Hence, when applying $\Omega(\Pi)$ to $f$, if $\sigma_S(T)=(v,q)$, it will be added to $X$. Thus, the first edge to be added to $X$ when applying $\Omega(\Pi)$ to $f$ will be in $E(T_n)$. 

We make a few observations.
\begin{itemize}
    \item \textit{Observation A:} For any $v\in U_k$, we know that if $e<_{\pi(T_{k-1})}e'$, $e$ and $e'$ both edges from $v$ to $T_{k-1}$, then $e<_{\pi(T_{k})}e'$, and $e,e'<_{\pi(T_{k})}e''=(p_k,v)$, if such an edge exists.
    \item \textit{Observation B:} When $v\in W_k$, we know the $f(v)$ edges which the map $\Phi_{\Pi}$ ``ignores''. That is, the set of smaller edges in condition $(2)$. Call this set $E_v$.
    \item \textit{Observation C:} The edge $(v,u)$ which will eventually connect $v$ to $T_k$ for some $k$ is determined as soon as $v\in W_k$.
\end{itemize} Elaborating on observation C, suppose $v\in W_k$ for $m\leq k\leq m+i$; i.e. $v$ is in $W_k$ for the first time when $k=m$, and is added to $T$ when $k=m+i$. Then by observation A, the order on the set $E_v\cup e_{m+i}$ is immutable for each of these $W_k$. In particular, the edges in $E_v$ are always smaller than $e_{m+i}$. Thus, only condition $(3)$ is not satisfied until $k=m+i$. At step $m+i$, the edge $e_{m+i}$ is smaller than any edge from $v$ to $T_{m+i}$ that is not in $E_v$. 

Assume by induction that thus far $T=(U,X)\subset T_n$. Then the next edge $e=\sigma_S(T)$ that is added to $X$ will be in $T_n$. Indeed, suppose $\sigma_S(T)=e$ is deleted. Then we know $e$ satisfies $(i),(ii)$, but $f(e_0)\geq 1$. This is true until $\sigma_S(T)=e'$ is added to $X$. The edge $e'$ is greater than exactly $f(e'_0)$ edges from $e_0$ to $T$ by observation A, and by observations B and C, we know that $e'$ must be in $E(T_n)$.
\end{proof}

\noindent Note that there may be several ways to define a map $\{\Pi_G\} \rightarrow \{\Sigma\}$. However, $\Omega$ was specifically defined so that it is injective. 

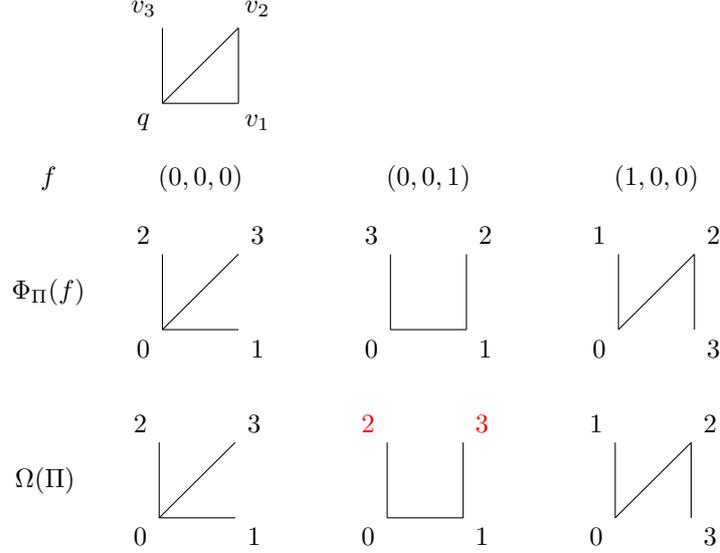
\begin{figure}[htb]
\centering
\begin{tikzpicture}

\draw (3,3)--(4,3);
\draw (2.75,2.75) node {$q$};
\draw (3,3)--(3,4);
\draw (4,3)--(4,4);
\draw (4.25,2.75) node {$v_1$};
\draw (2.75,4.25) node {$v_3$};
\draw (4,4)--(3,3);
\draw (4.25,4.25) node {$v_2$};

\draw (1.5, 2) node {$f$};

\draw (1.5,0.5) node {$\Phi_{\Pi}(f)$};

\draw (3.5,2) node {$(0,0,0)$};
\draw (3,0)--(4,0);
\draw (2.75,-.25) node {$0$};
\draw (3,0)--(3,1);
\draw (4.25,-.25) node {$1$};
\draw (2.75,1.25) node {$2$};
\draw (4,1)--(3,0);
\draw (4.25,1.25) node {$3$};

\draw (6.5,2) node {$(0,0,1)$};
\draw (6,0)--(7,0);
\draw (5.75,-.25) node {$0$};
\draw (6,0)--(6,1);
\draw (7.25,-.25) node {$1$};
\draw (7,0)--(7,1);
\draw (5.75,1.25) node {$3$};
\draw (7.25,1.25) node {$2$};

\draw (9.5,2) node {$(1,0,0)$};
\draw (8.75,-.25) node {$0$};
\draw (9,0)--(9,1);
\draw (8.75,1.25) node {$1$};
\draw (10,0)--(10,1);
\draw (10.25,1.25) node {$2$};
\draw (10,1)--(9,0);
\draw (10.25,-.25) node {$3$};

\end{tikzpicture}

\vspace{5mm}

\begin{tikzpicture}

\draw (1.5,0.5) node {$\Omega(\Pi)$};

\draw (3,0)--(4,0);
\draw (2.75,-.25) node {$0$};
\draw (3,0)--(3,1);
\draw (4.25,-.25) node {$1$};
\draw (2.75,1.25) node {$2$};
\draw (4,1)--(3,0);
\draw (4.25,1.25) node {$3$};

\draw (6,0)--(7,0);
\draw (5.75,-.25) node {$0$};
\draw (6,0)--(6,1);
\draw (7.25,-.25) node {$1$};
\draw (7,0)--(7,1);
\draw [color=red] (5.75,1.25) node {$2$};
\draw [color=red] (7.25,1.25) node {$3$};

\draw (8.75,-.25) node {$0$};
\draw (9,0)--(9,1);
\draw (8.75,1.25) node {$1$};
\draw (10,0)--(10,1);
\draw (10.25,1.25) node {$2$};
\draw (10,1)--(9,0);
\draw (10.25,-.25) node {$3$};

\end{tikzpicture}
\caption{The order in which vertices are added does not match between the maps.}
\label{difft}
\end{figure}

We observe that the order in which the vertices are added to $T_f$ according to $\Omega(\Pi)$ may not be the same as the order $\pi(T_f)$. An example is shown in Figure \ref{difft}. Let $\Pi_G$ be the proper set of tree orders in Example 3.2.2. The top row shows the global order on the vertices. The middle shows the $G$-parking functions and their images under $\Phi_{\Pi}$. The bottom row is the order in which the vertices are added when applying $\Omega(\Pi)$ to the corresponding functions. Nonetheless, the bijection between $\mathcal{P}_{G,q}$ and $\mathcal{T}_G$ is the same.
\subsection{Process Orders}
\label{cmy}
A bijection between G-parking functions and monomials of the Tutte polynomial that does not go through spanning trees was constructed by Chang et al \cite{CHANG2010231}. We describe this bijection and compare it to a tree growing sequence. Fix a total order $O_V:V\rightarrow \{0,\dots,n\}$ on the vertices of $G$. For each $f\in \mathcal{P}_{G,q}$, associate a \textit{process order} $\pi_f$ \cite{KOSTIC200873}. This is done recursively as follows:
\begin{enumerate}
   \item Let $\pi_f(0)=v_0=q$, and $V_0=V(G)-\{q\}$.
    \item Let $\pi_f(i)=min_{O_V}\{w\in V_{i-1}\:|\:0\leq f(w)< outdeg_{V_{i-1}}(w)\}$, where the vertices in $V-V_{i-1}$ have been processed. 
    \item Increase $i$ by one, and repeat step $2$ until all vertices have been processed; i.e. when $i=n-1$.
\end{enumerate}

One can get the process order for $f\in \mathcal{P}_{G,q}$ from a tree growing sequence. If $T$ does not span $S$, define $\Sigma_{O_V}$ by $\sigma_S(T)=(v,u)$, such that $v=min_{O_V}\{w\in V(S)-V(T)\}$, and $u$ is the smallest neighbor of $v$ in $T$. If $T$ spans $S$, then define $\sigma_S(T)$ to be the smallest edge according to the lexicographic order induced by $O_V$. Apply $\Sigma$ to any $f$. When the $i$-th edge is added to $X=E(T)$, identify $V_{i+1}$ with the vertices in $V-V(T)$. If $\sigma_S(T)=(v,u)$ is the edge added to $X$, then $\pi_f(|X|)=v$. Thus, $v$ will be added to $V(T)$ when at least $f(v)$ edges from $v$ to $T$ have been deleted and it is the minimum among all such candidate vertices, which is exactly statement $2$ above.

Denote $K=\{u\in V(G)\:|\:\pi^{-1}(v)\leq \pi^{-1}(u)\}$. This leads to the definition of a critical bridge vertex of ~$f$.\medskip
\\
\noindent \textbf{Definition 3.3.1:} \textit{A critical bridge vertex $v$ of the parking function $f$ with $\pi_f(i)=v$ is one for which $outdeg_K(v)=f(v)+1$ in G (criticality), and for every parking function $h$ satisfying: $g(\pi_h(j))=f(\pi_f(j))$, with $\pi_h(j)=\pi_f(j)$ for $j< i$; $h(v)\geq f(v)$; and $\pi_h(i)\geq_{O_V} v$; we have, in fact, that $\pi_h(i)=_{O_V}v$.}
\medskip
\newline
The last inequality says that there is no vertex strictly greater than $v$ according to $O_V$ which is processed at the same step as $v$ for some other $G$-parking function $h$. Let $cb_G(f)$ be the number of critical bridge vertices of $f$, and $w_G(f)=|E|-|V|+1-\sum_{v\in V-\{q\}} f(v)$.\medskip
\\
\textbf{Theorem 3.3 \cite{CHANG2010231}: }\textit{The Tutte polynomial of $G$ with fixed root $q$ can be expressed as the following closed formula: $$T(G;x,y)=\sum_{f\in \mathcal{P}_{G,q}} x^{cb_G(f)}y^{w_G(f)}$$}

\noindent Any tree growing sequence can be viewed as a way to write such a closed formula from the bijection $\rho: \mathcal{P}_{G,q}\rightarrow \mathcal{M}_G$. Simply say $T(G;x,y)=\sum_{f\in \mathcal{P}_{G,q}}\rho(f)$. However, any bijection comes with another bijection $\tau: \mathcal{P}_{G,q}\rightarrow \mathcal{T}_G$. We think that the above theorem secretly constructs a spanning tree and can be obtained via some tree growing sequence $\Sigma$. Specifically, the tree growing sequence $\Sigma_{O_V}$ defined above is the most likely candidate, and our conjecture has evidence through several calculations. However, we have not translated the constructions in \cite{CHANG2010231} to our language of tree growing sequences, and at this point we cannot verify the conjecture.

\section{Zonotopes}
\subsection{The Tutte Polynomial of a Zonotopal Tiling}
\label{zonotopes}

In the same spirit as the tree growing sequence algorithm, we describe a splitting algorithm which can be used to obtain the Tutte polynomial for a tiling of a zonotope. Let $M$ be an $n$-dimensional vector space, $N=M^{\vee}$, and $<,>$ the pairing of $N$ with $M$ (viewed as the standard inner product on $\mathbb{R}^n$).\medskip
\\
\textbf{Definition 4.1.1:} A \textbf{zonotope} is the image of a d-dimensional cube $Q_d=[0,1]^d$ under an affine projection. Equivalently, it is a Minkowski sum $$Z=\{a_1v_1+\dots+a_dv_d\,|\, 0\leq a_i\leq1,\,v_i\in M\cong \mathbb{R}^n\}.$$ We say that the set $X=\{v_1,\dots,v_d\}$ \textit{generates} $Z$.
\vspace{2mm}
%

\noindent We now discuss zonotopal tilings, where many of the details can be found in \cite{ProCon}, \cite{RichterGebertZiegler1993}, and ~\cite{ziegler95}.\medskip
\\
\textbf{Definition 4.1.2:} A \textbf{parallelotope} is a zonotope generated by vectors which form a basis of $M\cong \mathbb{R}^n$.\medskip
\\
\textbf{Definition 4.1.3:} A \textbf{cubical zonotopal tiling $\mathcal{Z}$} of $Z$ is a polyhedral complex comprised of a finite number of zonotopes $\{Z_i\}$ such that the maximal dimensional zonotopes - called \textbf{tiles} - are parallelotopes, and $\bigcup_i Z_i=Z$.
\vspace{2mm}

Let $Z$ be generated by $\{v_1,\dots,v_d\}$, and let $\{\mathcal{E}_j\}$ be the equivalence classes of edges of a tiling $\mathcal{Z}$, where an equivalence class is generated by the edges which are opposite a 2-dimensional face of $\mathcal{Z}$. Pick a representative $w_j$ of each $\mathcal{E}_j$. Each $w_j$ is parallel to a vector in the generating set $\{v_1,\dots,v_d\}$. As a result, we break each $v_i$ in to a finite number of vectors $w_1^i=k_1^iv_i,\:w_2^i=k_2^iv_i,\dots,\:w_l^i=k_l^iv_i$, such that $k_i>0$, $\sum k_i=1$, and $\sum_jw_j^i=v_i$. Then we can associate to $\mathcal{Z}$ the vector configuration $V_{\mathcal{Z}}$ containing the vectors $\{w_j^i\}$ for all $1\leq i\leq d$.\medskip
\\
\noindent \textbf{Definition 4.1.4:} Let $\mathcal{E}$ be an equivalence class of edges as above with representative $w$. A \textbf{zone} $B_w$ of a zonotope $Z$ is the set of tiles which contain an edge in $\mathcal{E}$. Two zones $B_w,\: B_{w'}$ are \textbf{parallel} if $w'$ is parallel to $w$.
\vspace{2mm}

\noindent Each zone $B_w$ has a \textit{positive side} $Z^{w,+}$ and \textit{negative side} $Z^{w,-}$ according to the direction of the vector $w$. An example is shown below in Figure 7.
\vspace{2mm}

\begin{figure}[htb]
    \centering
    \includegraphics[scale=0.3]{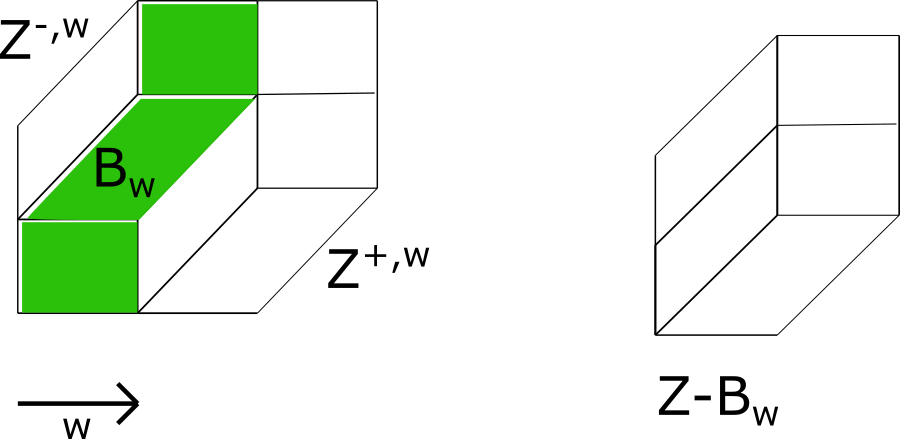}
    \caption{Shrinking.}
    \label{shrink}
\end{figure}

Fix a cubical tiling $\mathcal{Z}$ of $Z$. We compute a polynomial $T^{\ast}(\mathcal{Z};x,y)$ using a splitting algorithm which assigns a monomial to each tile of $\mathcal{Z}$, and $T^{\ast}(\mathcal{Z};x,y)$ is the sum of these monomials. We will be performing two operations - called \textit{shrinking} and $\textit{projection}$ in \cite{RichterGebertZiegler1993} - which will split the set of tiles into two disjoint sets at each step.\medskip
\\
\textbf{Definition 4.1.5:} Delete the zone $B_w$ and glue the positive and negative sides (see figure ~\ref{shrink}) of $B_w$ along $B_w\cap Z_w^+$ and $B_w\cap Z_w^-$. Denote the result of this operation $Z-B_w$ the \textbf{shrinking} of $Z$ with respect to $B_w$. Explicitly, since $w$ is parallel to $kv_i$ for some $0< k\leq 1$, we can write this zonotope as $$Z-B_w=\{a_1v_1+\dots+a_i(1-k)v_i+\dots a_nv_n\:|\:0\leq a_i\leq 1\}.$$
The tiling of $Z-B_w$ is as before. The associated vector configuration is $V_{\mathcal{Z}}-\{w\}$.\medskip
\\
\textbf{Definition 4.1.6:} Define $P_w:M\rightarrow M/(\mathbb{R}\cdot w)$. Let $Z|B_w=P_w(B_w)$ be the \textbf{projection} of the zone $B_w$. The tiles of $Z|B_w$ are $\{P_w(Z_i) \:| \: Z_i \text{ a tile of } B_w \}$. The associated vector configuration is $(V_{\mathcal{Z}}-\{w\})/(\mathbb{R}\cdot w)$.
\vspace{2mm}

Note that there is a description of $Z$ in terms of $Z-B_w$ and $Z|B_w$ in \cite{ProCon}. Our description is essentially the same, except we keep track of the tilings at each step. The decomposition
$$Z=(Z-B_w)\cup B_w$$
tells us that the set of tiles of $\mathcal{Z}$ splits into the tiles of $Z-B_w$ and tiles of $Z|B_w$.

Start with the tuple $(Z,\alpha,\beta)$, 
where initially $\alpha=0$, and $\beta=0$. The monomials will be $x^{\alpha}y^{\beta}$ where the exponent values will change according to the algorithm. Choose a belt $B_w$, and apply the shrinking and projection operations. This results in two tuples $(Z-B_w,\: 0,\: 0)$ and $(Z|B_w,\: \gamma,\: 0)$ associated to the resulting zonotopes, where $\gamma$ is the number of zones parallel to $B_w$. If $Z-B_w$ and $Z|B_w$ are zonotopes with the same tiling (i.e. they are equivalent zonotopes), we get a single tuple $(Z-B_w,\: 0,\:1)$. 
\begin{figure}[htb]
    \centering
    \includegraphics[scale=0.4]{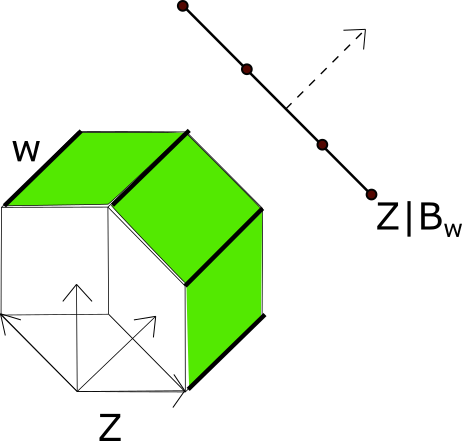}
    \caption{Projection with respect to $w$.}
    \label{union}
\end{figure}
Similar to the TGS algorithm, we repeat the operations for each new zonotope created. When the zonotope has been reduced to a collection of points with assigned tuples $(\bullet, \alpha, \beta)$, define $T^{\ast}(Z;x,y)$ to be the polynomial obtained by summing the monomials $x^{\alpha}y^{\beta}$. If we follow the path according to the splitting from each tile of $\mathcal{Z}$ to a point, we can associate a monomial to each tile. Thus, we can write the closed formula $$T^{\ast}(\mathcal{Z};x,y)=\sum_{\text{tiles of} \mathcal{Z}}x^{\alpha}y^{\beta}$$


Written in parallel to the deletion/restriction definition of the Tutte polynomial, and defining $T^{\ast}(\bullet;x,y)=~1$, the algorithm gives us the recursive formula:
\begin{equation}T^{\ast}(\mathcal{Z};x,y)=\begin{cases} 
      yT^{\ast}(\mathcal{Z}-B_w;x,y) & Z-B_w\cong Z|B_w\\
      x^{\gamma}T^{\ast}(\mathcal{Z}|B_w;x,y)+T^{\ast}(\mathcal{Z}-B_w;x,y) & \text{otherwise}
   \end{cases}
\end{equation}

\noindent \textbf{Observation:} The exponent $\gamma$ can be expressed in terms of vector configurations as\\ $\gamma~=|V_{\mathcal{Z}}|-|(V_{\mathcal{Z}}-\{w\})/\{w\}|$, the number of $0$-vectors resulting from the projection operation. We will `ignore' these $0$-vectors after projection, and can think of removing them from the configuration.\medskip
\\
\textbf{Example 4.1.1:} The zonotope generated by the vectors $v_1=(2,0), v_2=(0,1.5),\\v_3=(1,1)\in M\cong \mathbb{R}^2$ is a hexagon. Let $V_{\mathcal{Z}}=\{\frac{1}{2}v_1,\frac{1}{2}v_1,\frac{1}{2}v_2,\frac{1}{2}v_2, v_3\}$ be the vector configuration arising from the cubical zonotopal tiling $\mathcal{Z}$ shown below.

\begin{center}
\begin{tikzpicture}

\draw[->][very thick] (0,0)--(0,1.5);
\draw[->][very thick]  (0,0)--(1,1);
\draw[->][very thick]  (0,0)--(2,0);
\draw (0,1.5)--(1,2.5);
\draw (2,0)--(3,1);
\draw (3,1)--(3,2.5);
\draw (1,2.5)--(3,2.5);
\draw (1,0)--(1,0.75);
\draw (1,0.75)--(0,0.75);
\draw (1,0.75)--(2,1.75);
\draw (1,0)--(2,1);
\draw (2,1)--(3,1);
\draw (2,1.75)--(3,1.75);
\draw (0,0.75)--(1,1.75);
\draw (2,1)--(2,1.75);
\draw (1,1.75)--(1,2.5);
\draw (1,1.75)--(2,1.75);
\draw (2,1.75)--(2,2.5);

\end{tikzpicture}
\end{center}

\begin{figure}[ht]
\centering
    \includegraphics[scale=0.19]{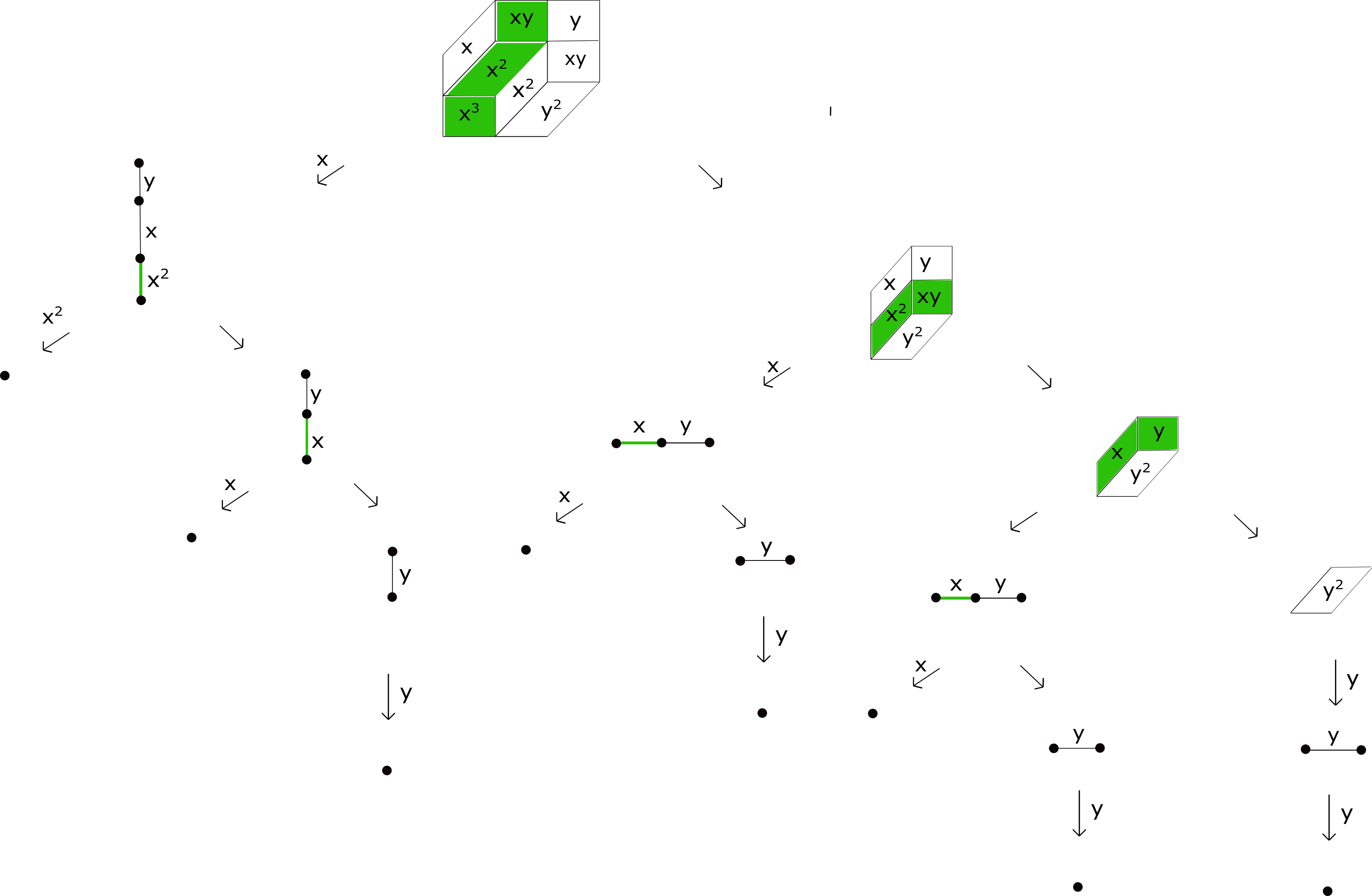}
    \vspace{3mm}
    
    \caption{$T^{\ast}(\mathcal{Z};x,y)=x^3+2x^2+x+2xy+y+y^2$}
    \label{zone}
\end{figure}

\noindent The splitting algorithm for $\mathcal{Z}$ is shown in Figure \ref{zone}, where in the first step the belt $B_w$ is chosen, where $w$ is the first $\frac{1}{2}v_1$ in the list. A colored zone means we are shrinking/projecting along that zone. A southwest arrow indicates projection, a southeast arrow indicates shrinking, and a south arrow represents when both are equivalent. Each intermediate zonotope $Z_k$ is tiled and the tiles labeled by the corresponding monomials of $T^{\ast}(\mathcal{Z}_k;x,y)$. The arrows are labeled according to where we multiply $T^{\ast}(\mathcal{Z}_k;x,y)$ by $x^{\gamma}$ or $y$ in the algorithm. The polynomial is $T^{\ast}(\mathcal{Z};x,y)=x^3+2x^2+x+2xy+y+y^2$, which is the Tutte polynomial for the graph $K_4-\{edge\}$.
\vspace{2mm}

\noindent We now recall a few notions related to matroids.\medskip
\\
\textbf{Definition 4.1.7:} The \textbf{rank function} of a matroid $M=(E, \mathcal{I})$ is 

$$r:2^E\rightarrow \mathbb{Z}_{>0}$$
$$r(A)=\max_{I\subseteq A, I\in \mathcal{I}}|I|$$

\noindent \textbf{Definition 4.1.8:} The \textbf{dual matroid $M^{\ast}$} is the pair $(E,\mathcal{I}^{\ast})$, where a set $J\in \mathcal{I}^{\ast}$ is independent in $M^{\ast}$ if and only if $E-J$ contains a basis of $M$. The rank function is  $r^{\ast}(A)=|A|-r(E)+r(E-A)$ is the dual rank function. 
\vspace{2mm}

\noindent \textbf{Definition 4.1.9:} The \textbf{Tutte polynomial} of a matroid is defined as 
\begin{equation}T(M;x,y)=\sum_{A\subseteq E}(x-1)^{r(E)-r(A)}(y-1)^{|A|-r(A)}\end{equation}
If $M^{\ast}$ is the matroid dual, then $T(M^{\ast};y,x)=T(M;x,y)$. Evaluating $T(M;2,2)$ gives the number of bases of $M$.
\vspace{2mm}

\noindent \textbf{Example 4.1.2:} If $G=(V,E)$ is a connected graph, we can define the \textit{cographical matroid} to be the matroid with ground set $E$ and bases $\mathcal{B}=\{\mathbf{b}=E-E(T)\,|\, $T$ \text{ a spanning tree}\}$. Hence, its rank is the genus $g=|E|-|V|+1$. More thorough expositions on matroids and their duals can be found in the original paper by Whitney \cite{WHIT35}, and lectures by Tutte \cite{Tutte1965LecturesOM}.

Let $w:~E(G)\rightarrow ~\mathbb{R}_{>0}$ be a function assigning length one to every edge of $G$, so that each edge can be identified to a unit interval. The zonotope $Z(G)$ is the projection of the cube $[0,1]^{|E(G)|}$ along the lattice of bonds (minimal cut-sets), which are the circuits of the cographical matroid. The dimension of $Z(G)$ is $g$.
Choose a tiling so that the representatives $w$ correspond to the edges of $G$. Then every tile corresponds to an element $\mathbf{b}\in \mathcal{B}$, the complement of a spanning tree. Hence, each tile corresponds to a unique spanning tree, and we get that any zone $B_w$ is the set of tiles associated to spanning trees which do not contain the edge corresponding to ~$w$.
\vspace{2mm}

\vspace*{\fill}

\vspace*{\fill}
\newpage

Let $\mathcal{W}$ be a matroid with ground set the list of vectors $W=\{w_1,\dots,w_d\}$ spanning the vector space $U=\mathbb{R}W$, and with independent sets $\mathcal{I}=\{X\, | \,W-X \text{ is a linearly independent set}\}$. Define $W_1=W-\{u\}$ and $W_2=W_1/u$. Then the Tutte polynomial of $W$ satisfies 
\begin{equation}T(W;x,y)=\begin{cases} 
      xT(W_1;x,y) & w \text{ is a coloop } (w=0)\\
      yT(W_2;x,y) & U=\mathbb{R}W_1\oplus \mathbb{R}\cdot w \text{ (w is a loop)}\\
      T(W_1;x,y)+T(W_2;x,y) & U=\mathbb{R}U_1,\:w\neq0
   \end{cases}
\end{equation}
The formula reduces to computing the Tutte polynomial of lists of vectors $V^{(i)}=V^{(i)}_0\sqcup V^{(i)}_1$, where $V^{(i)}_0$ is a list of $k$ linearly independent vectors and $V^{(i)}_0$ is a list of $h$ zero vectors; for such lists, $T(V^{(i)};x,y)=x^hy^k$. The bases of W are complements of subsets which form a basis for $U$. See \cite{MOCI09}, for example, for a treatment of the Tutte polynomial and a multiplicity polynomial for a vector configurations, as well as a discussion of how these polynomials give information about the associated zonotopes.\smallskip
\\
\textbf{Observation:} Given a zonotopal tiling $\mathcal{Z} $, the bases of the matroid $\mathcal{W}$ with ground set $V_{\mathcal{Z}}$ described above are in bijection with the tiles, since each tile has edges which form a basis for $U=\mathbb{R}V_{\mathcal{Z}}$. We will denote this matroid by $V_{\mathcal{Z}}^{\ast}$.\medskip
\\
\textbf{Theorem 4.1.1:} \textit{Fix a cubical zonotopal tiling $\mathcal{Z}$ of $Z$ with associated vector configuration $V_{\mathcal{Z}}$. Then $T^{\ast}(\mathcal{Z};x,y)$ is the Tutte polynomial $T(V_{\mathcal{Z}}^{\ast};x,y)$.}
\begin{proof}
Suppose we compute the Tutte polynomial $T(V^{\ast}_{\mathcal{Z}};x,y)$ and the polynomial
$T^{\ast}(\mathcal{Z};x,y)$ simultaneously, where the choice of $w$ at each step is a \emph{nonzero} vector. If any $0$-vectors are created, we choose to remove them immediately from the list. The algorithm for computing $T^{\ast}(\mathcal{Z};x,y)$ gives a bijection $$\{\text{Tiles of } Z\}\leftrightarrow \{\text{monomials}\}.$$
Hence, both polynomials 
have the same number of monomials. Moreover, the operations of deletion and restriction  applied to $V_{\mathcal{Z}}$ with respect to $w$ yield precisely the vector configurations associated to the tilings of $Z-B_w$ and $Z|B_w$, respectively. 

Inductively, let $Y$ be a zonotope with tiling $\mathcal{Y}$ that is created at some step of the algorithm with assoicated vector configuration $V_{\mathcal{Y}}$. Choose $w\neq 0 \in V_{\mathcal{Y}}$. We check that the recursion formulas for the polynomials are the same in all cases.

If $Y-B_w\cong Y|B_w$, then $T^{\ast}(\mathcal{Y};x,y)=yT^{\ast}(\mathcal{Y}-B_w;x,y)=yT^{\ast}(\mathcal{Y}|B_w;x,y)$. This occurs when $Y$ is a prism of height $w$, so the vector $w$ is a loop of $V^{\ast}_\mathcal{Y}$. Then $T(V^{\ast}_{\mathcal{Y}};x,y)=yT(V^{\ast}_{\mathcal{Y}}-\{w\};x,y)$. 

If we project $Y$ with respect to $w$, we multiply $T^{\ast}(Z|B_w;x,y)$ by $x^{\gamma}$, where $\gamma$ is the number of belts parallel to $B_w$. Recall that this represents throwing out all $0$-vectors created by projection. Thus, the integer $\gamma$ is the number of coloops in $V_{\mathcal{Y}}/\{w\}$, and subsequently contracting all of them gives $T(V^{\ast}_{\mathcal{Y}}/\{w\};x,y)=x^{\gamma}T(V^{\ast}_{\mathcal{Y}}/\{w,0,\dots,0\};x,y)$. 

If we shrink $Y$ with respect to $w$, and $w$ is not a coloop of $W$, then the tiles of $Y-B_w$ have the same monomials associated to them as in $Y$. 

Hence,
$T^{\ast}(\mathcal{Y};x,y)=x^{\gamma}T^{\ast}(\mathcal{Y}|B_w;x,y)+T^{\ast}(\mathcal{Y}-B_w;x,y)=x^{\gamma}T^(V^{\ast}_{\mathcal{Y}}/\{w,0,\dots,0\};x,y)+T(V_{\mathcal{Z}}^{\ast}-\{w\};x,y)$. This proves that $T^{\ast}(\mathcal{Z};x,y)=T(V_{\mathcal{Z}}^{\ast};x,y)$.
\end{proof}

\noindent \textbf{Remarks:} 
\begin{enumerate}
    \item If a zonotope $Z$ is a prism of height $w$, then it is a Minkowski sum $Z=Z'+w$, where $w$ orthogonal to $Z'$; thus, $Z-B_w\cong Z|B_w$. The converse is also true. 
    \item The bijection between tiles and monomials is dependent on the order in which we choose $e$. Indeed, if we have the zonotope $Z$ where $Z$ is a segment with two tiles $e_1,\:e_2$, then choosing $e_1$ first will assign $x$ to $e_1$ and $y$ to $e_2$. Hence, we can switch the order and get the other possible assignment.
\end{enumerate}


\subsection{Tiles of a Zonotope and $G$-parking Functions}
\label{end}

Let $G$ be a graph with length function which assigns each edge length one, and fix a labeling $\{q, v_1,\dots, v_n\}$ on the vertices of $G$, where we have chosen a root $q$. Let $Div(G)=\mathbb{Z}^{|V(G)|}$ be the group of $\mathbb{Z}$-linear combinations of vertices, written as $f=(a_0, a_1,\dots,a_n)$, with $a_0$ the coefficient of $q$, and $a_i$ the coefficient of $v_i$ for all other $i$. The \emph{degree} of a divisor is $\sum a_i$, and $Div^k(G)$ denotes the set of divisors of degree $k$. 

Variations of the \emph{chip-firing} game can be played on the vertices of a graph. If $v_i\in V$ and\\ $f=(a_0,a_1,\dots,a_n)~\in Div(G)$, then the chip-firing move $\sigma_i$ is defined by
\begin{equation}
    \sigma_i(f)(v_j)=\begin{cases} a_j-deg(v_j) & i=j\\
    a_j+n(v,v_j) & i\neq j
    \end{cases}
\end{equation}
where $n(v,v_i)$ is the number of edges between $v_j$ and $v_i$. We say that two divisors $f$ and $g$ are \textbf{linearly equivalent}, written $f\sim g$, if $g$ can be obtained from $f$ via a sequence of chip-firing moves. A \emph{principal divisor} is one linearly equivalent to $0$. Every chip-firing move is a sum of $\sigma_i$'s, so that one can view linear equivalence as being generated by the cuts $b_v$, where $b_v$ is the set of edges incident to $v$.\medskip
\\
\textbf{Definition 4.2.1:} A \textbf{$q$-reduced divisor} $f$ is a $G$-parking function when restricted to $V(G)-~\{q\}$, and in addition has $f(q)=-\sum_{v\in V(G)-\{q\}}f(v)$. In particular, we can view $\mathcal{P}_{G,q}\subset Div^0(G)$.\medskip
\\
\textbf{Definition 4.2.2:} The \textbf{Picard group} (of degree 0) is $Pic^0(G):=Div^0(G)/(f\sim 0)$.\medskip
\\
Thus, the Picard group measures the failure of degree 0 divisors to be principal.\medskip
\\
\textbf{Theorem 4.2.1, \cite{BN06}, \cite{MZ08}:} \textit{There exists a unique q-reduced representative in every linear equivalence class of $Div^0(G)$.\medskip
\\}
\textbf{Corollary 4.2.2:} \textit{Elements of $Pic^0(G)$ are in bijection with elements of $\mathcal{P}_{G,q}$.\medskip}

We now relate the set of $G$-parking functions to integer points in the cographical zonotope $Z(G)$. 
We start with a discussion for which \cite{BSMF1} and \cite{Biggs1999} are used as primary references. 

Fix an arbitrary orientation on the edges of $G$, and write an edge as an ordered tuple $e=(e_h,e_t)$. Let $C_0(G;\mathbb{R})\cong\mathbb{R}^{|V((G)|}$ and $C_1(G;\mathbb{R})\cong \mathbb{R}^{|E(G)|}$ be the vector spaces of finite $\mathbb{R}$-linear combinations of the vertices and edges of $G$, respectively, called the \emph{$0$-chains} and \emph{$1$-chains}.  We have that $Div(G)=C_0(G;\mathbb{Z})$. There is a standard inner product on $C_1(G;\mathbb{R})$ given by $<\sum a_e e,\sum b_ee>=\sum a_eb_e$.

Consider the map
\begin{center}
    \begin{tikzcd}
        C_1 \arrow[r, "d"]& C_0
    \end{tikzcd}
\end{center}
where $d(\sum_ea_ee)=\sum_ea_e(e_t-e_h)$ is the usual differential. Hence, $d(C_1(G;\mathbb{Z}))=Div^0(G)$, as the image is generated by $d(e)=e_t-e_h$. Denote the $1$-cycles by $Z_1\cong H_1(G;\mathbb{R})$. Note that $Z_1$ is isomorphic to $\mathbb{R}^g$, where $g=|E|-|V|+1$. Its orthogonal complement in $C_1(G;\mathbb{R})$ is generated by the cuts $b_v$. Let $\Lambda=Z_1\cap C_1(G;\mathbb{Z})\cong H_1(G;\mathbb{Z})$. We call $\Lambda$ the \emph{lattice of integral cycles}. 

Let $P$ be the orthogonal projection below.

\begin{center}
\begin{tikzcd}
    C_1(G;\mathbb{R}) \arrow[r, "P"] & Z_1\supset \Lambda
\end{tikzcd}
\end{center}

\noindent \textbf{Definition 4.2.3:} The \textbf{Jacobian} of $G$ is the finite group

$$J(G)=\frac{P(C_1(G;\mathbb{Z}))}{\Lambda}$$

We map $Div^0(G)$ into the real torus $Z_1/\Lambda\cong H_1(G;\mathbb{R})/H_1(G;\mathbb{Z})\cong \mathbb{R}^g/\mathbb{Z}^g$ as follows. Choose a path $p_i$ - viewed as an element of $C_1(G;\mathbb{Z})$ - from the root $q$ to each vertex $v_i$, and lift $f=\sum a_iv_i$ to $d^{-1}(f)=\sum_ia_ip_i\in C_1(G;\mathbb{Z})$. Then apply the orthogonal projection $P$, and take the image modulo $\Lambda$.

We get a map
$$A: Div^0(G)\rightarrow Z_1/\Lambda$$
$$f\mapsto P(d^{-1}(f))(mod\, \Lambda)$$

This map is a discrete analog of the Abel-Jacobi map originating from complex algebraic geometry. The image is the Jacobian $J(G)$. It is well-defined, as choosing another path $p_i'$ from $q$ to some $v_i$ and lifting $f$ will result in a shift by an element in $\Lambda$ (that is, $p_i-p_i'\in \Lambda$). Furthermore, if two divisors are linearly equivalent, they are sent to the same point. To see this, observe that linearly equivalent divisors differ by a cut. This leads to the proof of the Abel-Jacobi theorem.\medskip
\\
\textbf{Theorem 4.2.3 (Abel-Jacobi)\cite{BSMF1}, \cite{BN06}:} \textit{The Abel-Jacobi map induces a group isomorphism between $Pic^0(G)$ and $J(G)$.}
\vspace{2mm}

\noindent \textbf{Corollary 4.2.4: }\textit{The real torus $Z_1/\Lambda$ contains $|P_{G,q}|$ integral points.} 
\vspace{2mm}

Hence, a fundamental domain for $Z_1/\Lambda$ contains $|P_{G,q}|$ integral points. In fact, the closure of a fundamental domain can be identified with the cographical zonotope $Z(G)$, and we can tile $Z(G)$ as mentioned in example 4.1.2 (see \cite{AN14} for details). We also get the correspondence between spanning trees and tiles. Choose a generic fundamental domain $Z(G)$ in $Z_1$. The genericity means that integral points will not be vertices of the zonotope, but will lie in the interior of the parallelotopes of the tiling, so that there is precisely one integral point in each tile. Thus, we get a bijection between $G$-parking functions and tiles of $Z(G)$. 

We end with a question which ties together the bijections discussed in this paper. These bijections depend on several choices. There is the choice of a root $q$ for $G$, which determines the set $\mathcal{P}_{G,q}$. There are possibly several ways to tile the zonotope $Z(G)$ to obtain the bijection between spanning trees and tiles. We have the choice of fundamental domain. There is also the choice of the order in which the algorithm for computing $T^{\ast}(\mathcal{Z};x,y)$ is applied. Additionally, we have a choice of a tree growing sequence $\Sigma$ from which we get the bijective maps $\rho$ and $\tau$. In light of Theorem 3.1.1, the TGS may contain an underlying choice such as a total order on the edges. Thus, we state the following:
\medskip
\\
\textbf{Question:} \textit{Are there choices which are compatible in that they make the diagram below commute?\medskip}

\begin{center}
    \begin{tikzcd}
        &&\mathcal{P}_{G,q} \arrow[rd, leftrightarrow] \arrow[lld] \arrow[dd, leftrightarrow]\\
        \mathcal{M}_G  &&& \{\text{Tiles of }Z(G)\}\arrow[lll]\\
        && \mathcal{T}_G \arrow[ru, leftrightarrow] \arrow[llu, below left]
    \end{tikzcd}
    \label{conjecture}
\end{center}
\vspace{1mm}

The map $\mathcal{T}_G\rightarrow \mathcal{M}_G$ is given by a total edge order $O_E$ from which we read off internal and external activities. The maps from $\mathcal{P}_{G,q}$ to $\mathcal{M}_G$ and $\mathcal{T}_G$ are the bijections $\rho$ and $\tau$ arising from a tree growing sequence $\Sigma$. We have double-headed arrows for where there are known invertible algorithms.\medskip
\\
\noindent \textbf{Theorem:} \textit{The lower triangle is commutative.}

\begin{proof} Fix a tiling $\mathcal{Z}$ of $Z(G)$ and a correspondence between edges of $G$ and elements of $V_{\mathcal{Z}}$, which produces the correspondence between tiles and spanning trees. Compute $T^{\ast}(\mathcal{Z};x,y)$, but keep track of additional data. For every edge $e\in E(G)$, let $z(e)$ be the number of times $e$ is a coloop (parallel to the element $w$ chosen) or a loop. Let $e_i<e_j$ if $z(e_i)>z(e_j)$. If $z(e_i)=z(e_j)$, arbitrarily choose which is larger. The resulting total order $e_{i_1}< \dots< e_{i_m}$ induces a bijection from $\mathcal{T}_G$ to $\mathcal{M}_G$ via external/internal activities, and will match the bijection induced by the zonotope algorithm. In other words, we are determining how active an edge is through this count. As a general rule, the earlier an edge is chosen in a path, the less active it will be, and the higher it is in the total order.
\end{proof}

\newpage
\bibliographystyle{hplain.bst}
\bibliography{tgspaper1}

\begin{thebibliography}{10}

\bibitem{AN14}
Yang An, Matthew Baker, Greg Kuperberg, and Farbod Shokrieh.
\newblock Canonical representatives for divisor classes on tropical curves and
  the matrix-tree theorem.
\newblock {\em Forum of Mathematics, Sigma}, 2:e24, 2014, 0707.3168.

\bibitem{BSMF1}
Roland Bacher, Pierre~de La~Harpe, and Tatiana Nagnibeda.
\newblock The lattice of integral flows and the lattice of integral cuts on a
  finite graph.
\newblock {\em Bulletin de la Soci\'et\'e Math\'ematique de France},
  125(2):167--198, 1997.

\bibitem{BN06}
Matthew Baker and Serguei Norine.
\newblock Riemann-roch and abel-jacobi theory on a finite graph.
\newblock {\em Advances in Mathematics}, 215(2):766 -- 788, 2007, math/0608360.

\bibitem{BAKER2013164}
Matthew Baker and Farbod Shokrieh.
\newblock Chip-firing games, potential theory on graphs, and spanning trees.
\newblock {\em Journal of Combinatorial Theory, Series A}, 120(1):164 -- 182,
  2013.

\bibitem{bernardi:hal-00088479}
Olivier Bernardi.
\newblock {A characterization of the Tutte polynomial via combinatorial
  embeddings}.
\newblock {\em {Annals of Combinatorics}}, 12(2):139--153, 2008.
\newblock 14 pages.

\bibitem{Biggs1999}
N.L. Biggs.
\newblock Chip-firing and the critical group of a graph.
\newblock {\em Journal of Algebraic Combinatorics}, 9(1):25--45, Jan 1999.

\bibitem{BRYOX92}
Thomas Brylawski and James Oxley.
\newblock {\em The Tutte Polynomial and Its Applications}, pages 123--225.
\newblock Encyclopedia of Mathematics and its Applications. Cambridge
  University Press, 1992.

\bibitem{CHANG2010231}
Hungyung Chang, Jun Ma, and Yeong-Nan Yeh.
\newblock Tutte polynomials and g-parking functions.
\newblock {\em Advances in Applied Mathematics}, 44(3):231 -- 242, 2010.

\bibitem{CHEBIKIN2005}
Denis Chebikin and Pavlo Pylyavskyy.
\newblock A family of bijections between g-parking functions and spanning
  trees.
\newblock {\em Journal of Combinatorial Theory, Series A}, 110(1):31 -- 41,
  2005.

\bibitem{CORI200344}
Robert Cori and Yvan~Le Borgne.
\newblock The sand-pile model and tutte polynomials.
\newblock {\em Advances in Applied Mathematics}, 30(1):44 -- 52, 2003.

\bibitem{DHAR90}
Deepak Dhar.
\newblock Self-organized critical state of sandpile automaton models.
\newblock {\em Phys. Rev. Lett.}, 64:1613--1616, Apr 1990.

\bibitem{KOSTIC200873}
Dimitrije Kosti$\acute{c}$ and Catherine~H. Yan.
\newblock Multiparking functions, graph searching, and the tutte polynomial.
\newblock {\em Advances in Applied Mathematics}, 40(1):73 -- 97, 2008.

\bibitem{MZ08}
G.~Mikhalkin and I.~Zharkov.
\newblock Tropical curves, their jacobians and theta functions.
\newblock {\em Contemporary Mathematics}, 465:\:203--230, 2008.

\bibitem{MOCI09}
Luca Moci.
\newblock A tutte polynomial for toric arrangements.
\newblock {\em Transactions of the American Mathematical Society}, 364, 11
  2009.

\bibitem{ProCon}
Claudio Procesi and Corrado De~Concini.
\newblock {\em Topics in Hyperplane Arrangements, Polytopes and Box-Splines}.
\newblock Springer-Verlag New York, 2010.

\bibitem{RichterGebertZiegler1993}
J{\"u}rgen Richter-Gebert and G{\"u}nter~M. Ziegler.
\newblock Zonotopal tilings and the bohne-dress theorem.
\newblock Technical Report SC-93-25, ZIB, Takustr. 7, 14195 Berlin, 1993.

\bibitem{tutte_1954}
W.~T. Tutte.
\newblock A contribution to the theory of chromatic polynomials.
\newblock {\em Canadian Journal of Mathematics}, 6:80--91, 1954.

\bibitem{Tutte1965LecturesOM}
William~T. Tutte.
\newblock Lectures on matroids.
\newblock {\em J. Res. Natl. Bur. Stand., Sec. B: Math. and Math. Phys.},
  69B(1-2):1--47, 1965.

\bibitem{WHIT35}
Hassler Whitney.
\newblock On the abstract properties of linear dependence.
\newblock {\em American Journal of Mathematics}, 57(3):509--533, 1935.

\bibitem{ziegler95}
Gunter~M. Ziegler.
\newblock {\em Lectures on Polytopes}, volume 152 of {\em Graduate Texts in
  Mathematics}.
\newblock Springer-Verlag New York, 1 edition, 1995.

\end{thebibliography}

\end{document}